\documentclass[11pt]{article}
\usepackage[a4paper,margin=2.5cm]{geometry}
\usepackage[T1]{fontenc}
\usepackage{lmodern,microtype}
\usepackage{amsmath,amssymb,amsthm}
\usepackage{booktabs}
\usepackage{graphicx}
\usepackage{xcolor}
\usepackage{tikz,pgfplots}
\usetikzlibrary{arrows.meta,positioning,calc,fit}
\pgfplotsset{compat=1.18}
\usepackage[font=small,labelfont=bf]{caption}
\usepackage{subcaption}
\usepackage[hidelinks]{hyperref}

\newtheorem{lemma}{Lemma}
\theoremstyle{remark}

\newcommand{\conv}{\operatorname{conv}}
\newcommand{\area}{\operatorname{area}}
\newcommand{\amin}{a_{\min}}
\newcommand{\R}{\mathbb{R}}

\newcommand{\StateDate}{28 November 2026}
\newcommand{\NOurs}{218}
\newcommand{\NCells}{294}
\newcommand{\NOursBig}{181}
\newcommand{\NBig}{192}
\newcommand{\NOursBigSquare}{53}
\newcommand{\CodeURL}{\url{https://github.com/spiralulam/heilbronn-pipeline}}

\definecolor{mplC0}{HTML}{1f77b4}
\definecolor{mplC1}{HTML}{ff7f0e}
\colorlet{accent}{mplC0}
\colorlet{accentdark}{mplC0!75!black}
\colorlet{hot}{mplC1}
\tikzset{
  cfgtri/.style={fill=accent,fill opacity=0.16,draw=accent,draw opacity=0.55,line width=0.25pt},
  cfgbox/.style={line width=0.8pt,black!80},
  cfgpt/.style={circle,fill=accentdark,inner sep=0pt,minimum size=3pt},
  cfghot/.style={circle,fill=hot,draw=white,line width=0.3pt,inner sep=0pt,minimum size=5pt},
  stripfill/.style={fill=accent,fill opacity=0.12},
  orb0/.style={circle,fill=accentdark,inner sep=0pt,minimum size=4.6pt},
  orb1/.style={circle,fill=hot,inner sep=0pt,minimum size=4.6pt},
  orb2/.style={circle,fill=black!70,inner sep=0pt,minimum size=4.6pt},
  orb3/.style={circle,draw=accentdark,fill=white,line width=0.7pt,inner sep=0pt,minimum size=4.4pt},
  stage/.style={rounded corners=3pt,draw=accent,fill=accent!6,line width=0.8pt,font=\small,align=center,
                minimum height=1.35cm,text width=2.5cm,execute at begin node={\hyphenpenalty10000\relax}},
  hotstage/.style={stage,draw=hot!85!black,fill=hot!8},
  flow/.style={-{Latex[length=2.2mm]},line width=0.9pt,draw=black!55},
}

\title{A record-winning optimization algorithm for the\\ Heilbronn triangle problem with AI assistance}
\author{Nathan Sudermann-Merx\thanks{Institute for Computer Science, DHBW Mannheim. Email: \href{mailto:nathan.sudermann-merx@dhbw.de}{nathan.sudermann-merx@dhbw.de}}
\and Tej Stead\thanks{ElevenLabs. Email: \href{mailto:tejstead@tejstead.com}{tejstead@tejstead.com}}}
\date{\StateDate}

\begin{document}
\maketitle

\begin{abstract}
The Heilbronn triangle problem asks for $n$ points in a convex region that maximize the smallest area of a
triangle formed by three of them. We describe a heuristic optimization pipeline for the three variants of the public
record tables (unit square, triangle, and the convex variant, in which the region is the convex hull of the points).
Its ingredients are a one-point relocation by linear programming, a min-norm ascent direction from nonsmooth analysis, stress-weighted ruin and recreate, symmetric
starts, moves between neighbouring $n$ and between the three variants, a GPU stage that screens tens of thousands of
such starts per minute, and a final polish consisting of trust-region sequential linear programming with a lossless
reduction of the constraint set, a high-precision solution of the system of smallest triangles, and verification in rational
arithmetic. As of \StateDate{} the pipeline holds the records in \NOurs{} of the \NCells{} cells with $n \le 100$,
including \NOursBig{} of the \NBig{} cells with $37 \le n \le 100$. The pipeline was developed, implemented and operated almost
entirely by an AI coding assistant; we release a compact version of the code.
\end{abstract}

\section{Introduction}

For a convex region $K \subset \R^2$ and a set $P$ of $n \ge 3$ points in $K$ let $\amin(P)$ be the smallest area of a
triangle with vertices in $P$. The Heilbronn problem asks for
\[
  \Delta_n(K) \;=\; \max_{P \subset K,\ |P| = n} \; \frac{\amin(P)}{\area(K)} .
\]
Asymptotically $\Delta_n$ lies between $c \log n / n^2$ \cite{KPS82} and $n^{-8/7-1/2000}$ \cite{CPZ23}, following
Roth's first non-trivial upper bound \cite{Roth51}. Exact values are known only for small $n$: mixed-integer nonlinear
optimization solves the square to global optimality (within numerical tolerances) for $n \le 10$ and reconstructs
exact coordinates \cite{SM26a}, and certifies
optimality in the triangle for $n \le 8$ up to a relative gap of $10^{-4}$ \cite{SM26b}. Beyond that the best known
configurations come from computation \cite{Goldberg72,CY02,AlphaEvolve}. The public record tables of the second author \cite{HSite}, a much extended
successor of Erich Friedman's page \cite{Friedman}, list $n = 3, \dots, 100$ for three variants: the unit square, the triangle, and the
\emph{convex variant}, where the region is the convex hull of the points themselves and the value is
$\amin(P)/\area(\conv P)$ (the points need not be in convex position). Every entry is given by exact decimal coordinates and verified in rational arithmetic.

\begin{figure}[t]
  \centering
  \includegraphics[width=\linewidth]{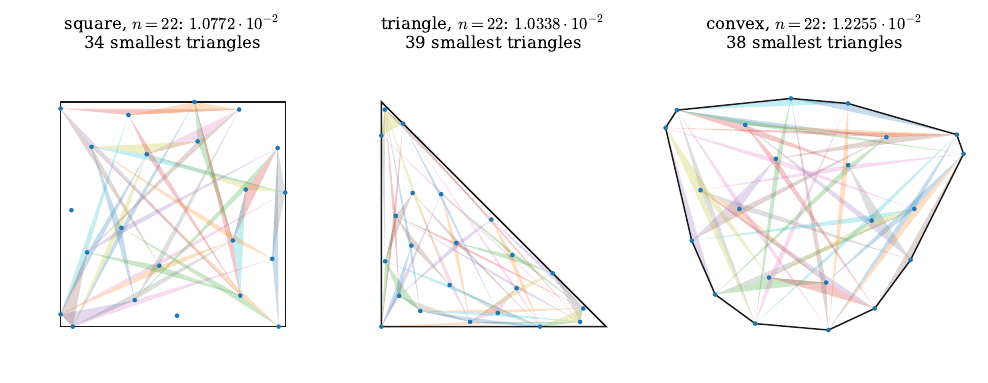}
  \caption{Our record configurations for $n = 22$ in the three variants, with every triangle of minimal area shaded,
  one colour per triangle (relative tolerance $10^{-6}$). Values are $\amin/\area$.}
  \label{fig:gallery}
\end{figure}

This paper describes the optimization algorithm behind \NOurs{} of the \NCells{} current records (\StateDate), among
them \NOursBig{} of the \NBig{} cells with $37 \le n \le 100$, including all $64$ triangle and all $64$ convex cells in
that range.
Figure~\ref{fig:gallery} shows three of them. We do not claim optimality for any configuration. The algorithm is a
heuristic built from simple, mostly classical parts; what makes it work is how the parts are combined.

\paragraph{Contributions and AI assistance.} The first author contributed some of the main ideas and set the rules of the work
(resource limits, verification, crediting). The second author built the public record tables as a successor of Erich
Friedman's page and extended them far beyond it, with an exact verifier for every entry, and with these tables
stimulated the community's interest in the problem. The second author was also involved in the search for the best
configurations from the beginning, provided the first configurations for many of the larger instances, and wrote the
first search engines, which served as the starting point for the later versions. Everything else, including the
remaining ideas, the implementation of the later versions, the experiments, the operation of the search around the
clock, the submission of results and the first draft of this paper, was done by an AI coding assistant (Claude,
Anthropic, through Claude Code) following the first author's instructions.

\section{Overview}\label{sec:overview}

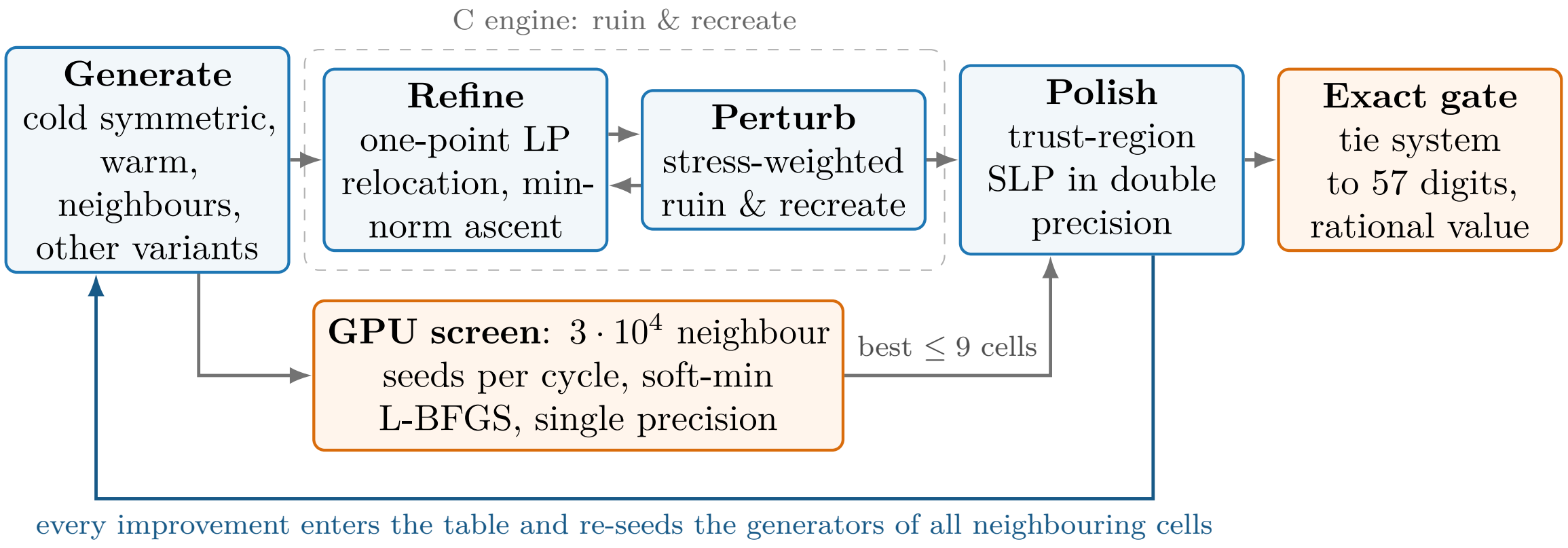
\begin{figure}[t]
  \centering
  \begin{tikzpicture}
    \node[stage] (g) at (0,0) {\textbf{Generate}\\cold symmetric, warm, neighbours, other variants};
    \node[stage] (r) at (3.1,0) {\textbf{Refine}\\one-point LP relocation, min-norm ascent};
    \node[stage] (p) at (6.2,0) {\textbf{Perturb}\\stress-weighted ruin \& recreate};
    \node[stage] (h) at (9.3,0) {\textbf{Polish}\\trust-region SLP in double precision};
    \node[hotstage] (e) at (12.4,0) {\textbf{Exact gate}\\tie system to 57 digits, rational value};
    \node[hotstage,text width=4.9cm,minimum height=1.0cm] (gpu) at (4.2,-2.1)
      {\textbf{GPU screen}: $3\cdot 10^4$ neighbour seeds per cycle, soft-min L-BFGS, single precision};
    \node[draw=black!30,dashed,rounded corners=4pt,fit=(r)(p),inner sep=5pt] (ce) {};
    \node[font=\scriptsize,text=black!60,anchor=south] at (ce.north) {C engine: ruin \& recreate};
    \draw[flow] (g) -- (r);
    \draw[flow] ([yshift=2.5mm]r.east) -- ([yshift=2.5mm]p.west);
    \draw[flow] ([yshift=-2.5mm]p.west) -- ([yshift=-2.5mm]r.east);
    \draw[flow] (p) -- (h);
    \draw[flow] (h) -- (e);
    \draw[flow] ([xshift=5mm]g.south) |- (gpu.west);
    \coordinate (hin) at ([xshift=-5mm]h.south); \coordinate (hout) at ([xshift=5mm]h.south);
    \draw[flow] (gpu.east) -- node[above,font=\scriptsize,text=black!70] {best $\le 9$ cells} (gpu.east -| hin) -- (hin);
    \draw[flow,draw=accentdark] (hout) -- (hout |- 0,-3.3) -- node[below,font=\scriptsize,text=accentdark]
      {every improvement enters the table and re-seeds the generators of all neighbouring cells} (-0.5,-3.3) -- ([xshift=-5mm]g.south);
  \end{tikzpicture}
  \caption{The pipeline. A polished result that beats the stored value by a relative $10^{-9}$ enters a table of best
  configurations, one per cell (variant, $n$); every generator reads its starting points from this table, so an
  improvement in one cell immediately changes the starts of its neighbours. The exact gate (Section~\ref{sec:tie}) is
  applied when a configuration is exported for submission.}
  \label{fig:pipeline}
\end{figure}

Figure~\ref{fig:pipeline} shows the pipeline. Four ideas carry it.

\begin{figure}[t]
  \centering
  \begin{tikzpicture}[
      cell/.style={rounded corners=2pt,draw=black!40,minimum width=1.35cm,minimum height=0.62cm,font=\small,fill=white},
      hotc/.style={cell,draw=accent,fill=accent!12,line width=1pt},
      ar/.style={-{Latex[length=1.8mm]},line width=0.9pt,draw=accentdark}]
    \foreach \r/\v in {0/square,1/triangle,2/convex}{
      \node[font=\small,anchor=east] at (-0.85,-\r) {\v};
      \foreach \c/\t in {0/{$n-2$},1.8/{$n-1$},3.6/{$n$},5.4/{$n+1$},7.2/{$n+2$}}{\node[cell] at (\c,-\r) {\t};}}
    \node[hotc] (m) at (3.6,-1) {$n$};
    \draw[ar] (m.north) -- (3.6,-0.31);
    \draw[ar] (m.south) -- (3.6,-1.69);
    \draw[ar] (m.west) -- (2.48,-1);
    \draw[ar] (m.east) -- (4.72,-1);
    \node[font=\scriptsize,text=black!60,align=left,anchor=west] at (8.3,-0.5) {vertical: affine maps (Lemma~\ref{lem:affine})};
    \node[font=\scriptsize,text=black!60,align=left,anchor=west] at (8.3,-1.5) {horizontal: delete or insert points};
  \end{tikzpicture}
  \caption{The table as one object: an improved cell (here triangle, $n$) provides starts for its neighbours along $n$
  (also at distance two) and in the other variants.}
  \label{fig:matrix}
\end{figure}
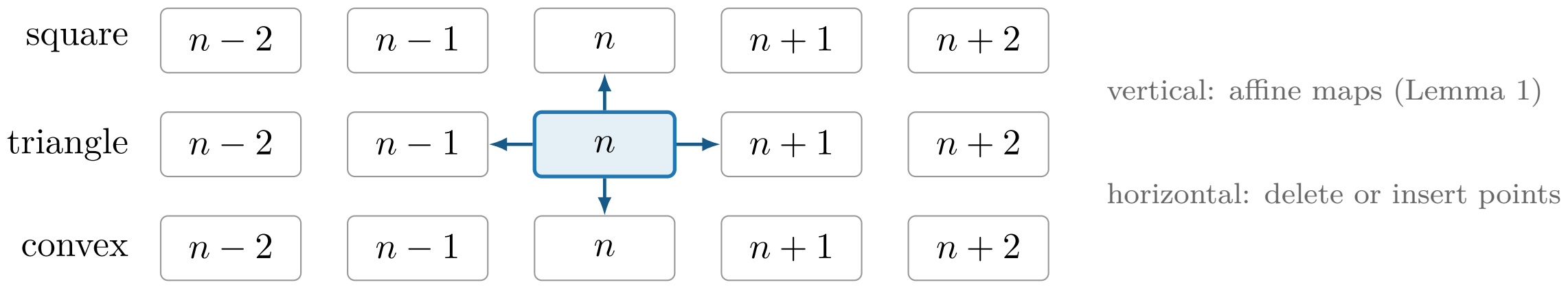

\paragraph{1. One table instead of 294 problems.} A good configuration with $n$ points carries most of the structure
of good configurations with $n \pm 1$ points: delete a well-chosen point, or insert one greedily, and a local
optimizer often reaches a better configuration than any search from scratch. Likewise a configuration in one variant
is a valid start in the others (Section~\ref{sec:moves}). We therefore treat the $3 \times 98$ cells as one connected
object (Figure~\ref{fig:matrix}); the generators read their starts from the current table and every improvement
radiates to the neighbouring cells.

\paragraph{2. The true objective at the end.} The objective is a minimum of $\binom{n}{3}$ quadratic functions, and at
a good configuration roughly as many of them are active as there are free parameters (Section~\ref{sec:tie}). Smooth
surrogates are useful for exploring, but the final steps work on the nonsmooth objective itself: a single point is
relocated to its best position in a face of a line arrangement by a three-variable LP, the polish solves LPs whose
constraint set is reduced without loss, and the final system of equal smallest triangles is solved to 57 digits before
the value is recomputed in rational arithmetic.

\paragraph{3. Screen in bulk, decide carefully.} Most starting configurations are useless, and the SLP polish is the
bottleneck. A GPU completes and soft-polishes tens of thousands of starts per cycle in single precision; only the best
start of a few promising cells reaches the SLP polish on the CPU (Section~\ref{sec:gpu}).

\paragraph{4. Aimed perturbation and symmetry.} Escapes from local maxima remove the points that carry the most
nearly smallest triangles (stress-weighted ruin and recreate), and cold starts are searched in the orbit space of a
point group, where good configurations are not of measure zero (Section~\ref{sec:global}).

\section{Setting}\label{sec:setting}

A configuration is $P = (p_1, \dots, p_n) \in \R^{2n}$, and $T = \binom{[n]}{3}$ is the set of triangles. For
$t = (i,j,k)$ let
\[
  A_t(P) \;=\; \tfrac12 \bigl( (x_j - x_i)(y_k - y_i) - (x_k - x_i)(y_j - y_i) \bigr),
  \qquad a_t = |A_t|, \qquad \amin = \min_{t \in T} a_t .
\]
$A_t$ is a quadratic polynomial, affine in each point separately, and $\sigma_t = \operatorname{sign} A_t$ denotes the
orientation of $t$ at the current configuration, so that $a_t = \sigma_t A_t$ is smooth nearby. The value of $P$ in a
region $K$ is $v(P) = \amin(P)/\area(K)$: $v = \amin$ for the unit square $[0,1]^2$, $v = 2\amin$ for the triangle
with vertices $(0,0),(1,0),(0,1)$, and $v = \amin/H$ with $H(P) = \area(\conv P)$ in the convex variant, where
$K = \conv P$. We write $v_t = a_t / \area(K)$, so that $v = \min_t v_t$.

\begin{lemma}\label{lem:affine}
Let $\phi(x) = Mx + b$ with $\det M \neq 0$. Then $A_t(\phi(P)) = \det M \cdot A_t(P)$ for all $t$ and
$\area(\phi(K)) = |\det M|\,\area(K)$, so every $v_t$ is invariant when $P$ and $K$ are mapped together. Moreover,
$\amin(P)/H(P) \ge \amin(P)/\area(K)$ whenever $P \subset K$.
\end{lemma}

The first part follows from the determinant formula, the second from $\conv P \subseteq K$. Hence every triangle is
equivalent to the reference triangle, every parallelogram to the square, the convex value is affine invariant, and a
square or triangle configuration is a convex configuration of at least the same value.

In the square and the triangle each point has a \emph{type}: interior (two free coordinates), on a side (one
coordinate along the side) or at a corner (none). Types are detected with tolerance $10^{-7}$ ($10^{-9}$ in the tie
refinement). The soft polish, the relocation and the tie refinement keep them fixed, which removes the container
constraints of boundary points; this matters because good configurations have many points on the boundary
(Figure~\ref{fig:gallery}). The SLP polish instead treats all coordinates as free and imposes the container as bounds,
so that points can reach or leave the boundary. In the convex variant all points are free. Of the six affine gauge
directions, translation and scale are normalized after every step (centroid $0$, $H = 1$); the other three do not
change the value, and only the tie refinement removes them (Section~\ref{sec:tie}).

\section{Local optimization}\label{sec:local}

\subsection{One point, one LP}\label{sec:reloc}

Fix all points except $p = (x,y)$. For a pair $j,k$ of the other points, $2A_{(p,j,k)} = c_{jk} + \alpha_{jk} x +
\beta_{jk} y$ is affine in $(x,y)$, so for a target value $z$ the condition $|A_{(p,j,k)}| \ge z$ forbids an open
strip of width $4z/\|p_j - p_k\|$ around the line through $p_j$ and $p_k$ (Figure~\ref{fig:reloc}). These lines cut
$K$ into convex polygons, the \emph{faces} of the arrangement, and inside a face all orientations
$\sigma_{jk} = \operatorname{sign} A_{(p,j,k)}$ are constant. The best position of $p$ in its current face therefore
solves the LP in three variables
\begin{equation}\label{eq:reloc}
  \max_{x,y,z}\ z \quad \text{s.t.} \quad \sigma_{jk}\,(c_{jk} + \alpha_{jk} x + \beta_{jk} y) \ge 2z \ \ \forall j<k,
  \qquad z \le \min_{t \not\ni p} a_t, \qquad (x,y) \in K .
\end{equation}
It is solved by Seidel's randomized algorithm \cite{Seidel91} in expected time linear in the number of constraints.
Since the current position is feasible, \eqref{eq:reloc} never decreases $\amin$, and a sweep over all points is
monotone; points on a side are moved along the side (an LP in two variables).\footnote{\emph{Multi-face probing} also
solves \eqref{eq:reloc} in the faces of about $40$ random positions and moves $p$ to the best of them, a
large-neighbourhood search over the placements of one point. It served as the re-insertion move in the first author's
earlier engine; the released engine uses greedy insertion (Section~\ref{sec:rr}) and the plain sweep.}

\subsection{Leaving a relocation-stationary point}\label{sec:ascent}

\begin{figure}[t]
  \centering
  \begin{subfigure}[b]{0.38\linewidth}\centering
    \begin{tikzpicture}[scale=4.6]
\begin{scope}
\clip (0,0) rectangle (1,1);
\fill[stripfill] (2.2272,-1.8389) -- (-1.7744,2.6318) -- (-2.2272,2.2265) -- (1.7744,-2.2442) -- cycle;
\fill[stripfill] (2.9882,0.4712) -- (-2.9986,0.0737) -- (-2.9882,-0.0836) -- (2.9986,0.3138) -- cycle;
\fill[stripfill] (0.5373,3.1468) -- (-0.7246,-2.7190) -- (-0.5373,-2.7593) -- (0.7246,3.1065) -- cycle;
\fill[stripfill] (-0.1249,3.1938) -- (-0.1249,-2.8062) -- (0.1249,-2.8062) -- (0.1249,3.1938) -- cycle;
\fill[stripfill] (2.9283,-0.4664) -- (-2.8759,1.0538) -- (-2.9283,0.8539) -- (2.8759,-0.6663) -- cycle;
\fill[stripfill] (2.0723,2.3641) -- (-2.1703,-1.8785) -- (-2.0723,-1.9766) -- (2.1703,2.2661) -- cycle;
\fill[stripfill] (2.4994,1.8544) -- (-2.5708,-1.3540) -- (-2.4994,-1.4668) -- (2.5708,1.7415) -- cycle;
\fill[stripfill] (3.0077,0.9877) -- (-2.7155,-0.8137) -- (-2.6608,-0.9877) -- (3.0624,0.8137) -- cycle;
\fill[stripfill] (0.0944,3.0000) -- (0.0944,-3.0000) -- (0.2525,-3.0000) -- (0.2525,3.0000) -- cycle;
\fill[stripfill] (-0.5342,2.9168) -- (0.6980,-2.9553) -- (0.8811,-2.9168) -- (-0.3511,2.9553) -- cycle;
\fill[stripfill] (3.1734,0.1395) -- (-2.8266,0.1395) -- (-2.8266,-0.1395) -- (3.1734,-0.1395) -- cycle;
\fill[stripfill] (1.7212,2.5708) -- (-1.4871,-2.4994) -- (-1.3743,-2.5708) -- (1.8340,2.4994) -- cycle;
\fill[stripfill] (2.2470,2.1691) -- (-1.9957,-2.0735) -- (-1.9001,-2.1691) -- (2.3426,2.0735) -- cycle;
\fill[stripfill] (-1.2829,2.2079) -- (3.1878,-1.7937) -- (3.2829,-1.6876) -- (-1.1878,2.3141) -- cycle;
\fill[stripfill] (-1.6443,1.6788) -- (3.5765,-1.2782) -- (3.6443,-1.1585) -- (-1.5765,1.7985) -- cycle;
\fill[stripfill] (-0.9695,-2.0130) -- (3.2732,2.2296) -- (2.9695,2.5334) -- (-1.2732,-1.7093) -- cycle;
\fill[stripfill] (0.1399,3.1361) -- (1.6602,-2.6681) -- (1.8601,-2.6158) -- (0.3398,3.1885) -- cycle;
\fill[stripfill] (0.8605,3.2602) -- (0.8605,-2.7398) -- (1.1395,-2.7398) -- (1.1395,3.2602) -- cycle;
\fill[stripfill] (-1.7201,-1.3491) -- (2.5226,2.8935) -- (2.0670,3.3491) -- (-2.1757,-0.8935) -- cycle;
\fill[stripfill] (1.7118,-1.5765) -- (-1.2452,3.6443) -- (-1.3649,3.5765) -- (1.5921,-1.6443) -- cycle;
\fill[stripfill] (3.1734,1.1249) -- (-2.8266,1.1249) -- (-2.8266,0.8751) -- (3.1734,0.8751) -- cycle;
\fill[stripfill] (3.1287,0.4755) -- (-2.7434,1.7077) -- (-2.7818,1.5245) -- (3.0903,0.2923) -- cycle;
\fill[stripfill] (2.0539,-1.3613) -- (-1.9477,3.1094) -- (-2.0539,3.0144) -- (1.9477,-1.4563) -- cycle;
\fill[stripfill] (2.9127,1.5512) -- (-2.9531,0.2893) -- (-2.9127,0.1019) -- (2.9531,1.3638) -- cycle;
\fill[stripfill] (3.0000,0.9056) -- (-3.0000,0.9056) -- (-3.0000,0.7475) -- (3.0000,0.7475) -- cycle;
\fill[stripfill] (0.8599,2.9986) -- (0.4624,-2.9882) -- (0.6198,-2.9986) -- (1.0172,2.9882) -- cycle;
\fill[stripfill] (1.5536,2.8890) -- (-0.2479,-2.8342) -- (-0.0739,-2.8890) -- (1.7275,2.8342) -- cycle;
\fill[stripfill] (3.2442,-0.7744) -- (-1.2265,3.2272) -- (-1.6318,2.7744) -- (2.8389,-1.2272) -- cycle;
\end{scope}
\draw[cfgbox] (0,0) rectangle (1,1);
\node[cfgpt] at (0.0000,0.1938) {};
\node[cfgpt] at (0.1734,0.0000) {};
\node[cfgpt] at (1.0000,0.2602) {};
\node[cfgpt] at (0.1734,1.0000) {};
\node[cfgpt] at (0.0000,0.8266) {};
\node[cfgpt] at (0.7398,0.0000) {};
\node[cfgpt] at (0.8062,1.0000) {};
\node[cfgpt] at (1.0000,0.8266) {};
\node[cfghot] at (0.6531,0.3469) {};\input{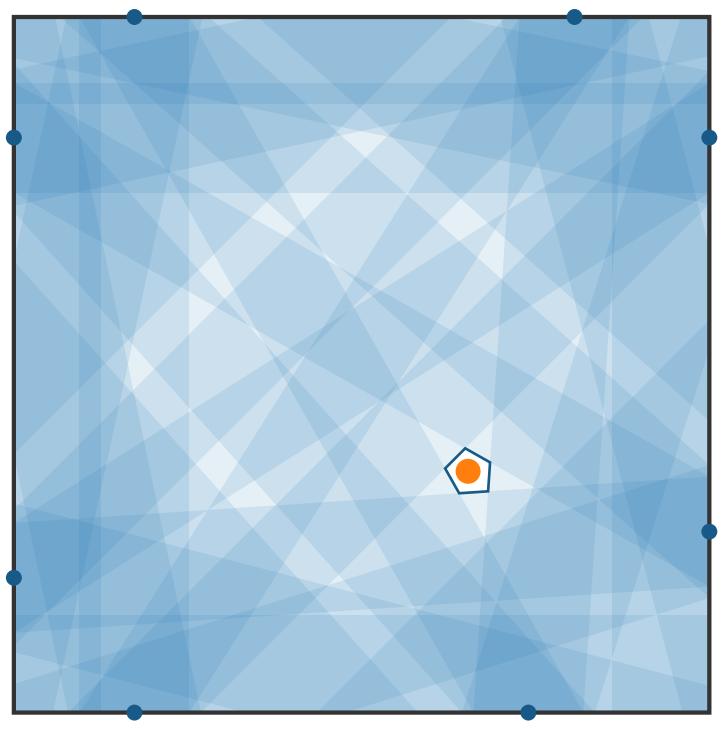}
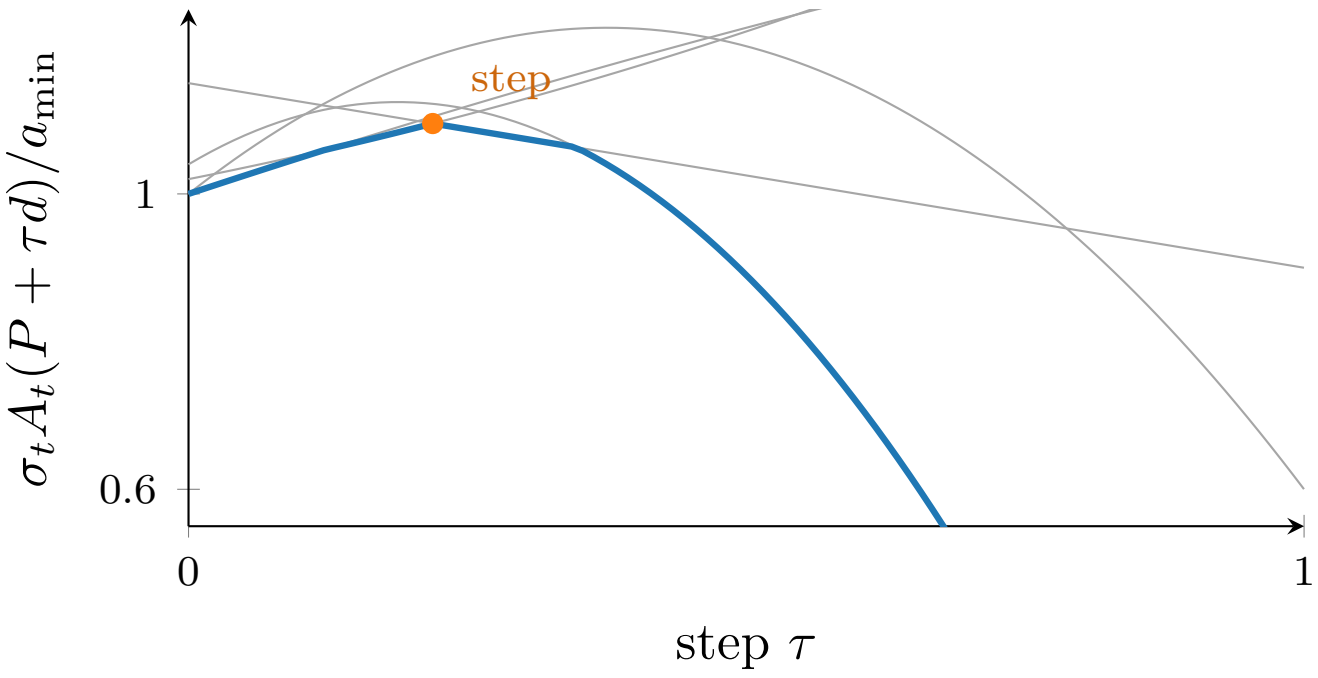
      \node[cfghot] at (0.6531,0.3469) {};\end{tikzpicture}
    \caption{Forbidden strips for one point (orange) of a configuration with $n=9$, for $z$ slightly below the current
    minimum. The outlined pentagon is the admissible part of its face, the feasible set of \eqref{eq:reloc}.}
  \end{subfigure}\hfill
  \begin{subfigure}[b]{0.56\linewidth}\centering
    \begin{tikzpicture}
      \begin{axis}[width=\linewidth,height=5.0cm,xmin=0,xmax=1,ymin=0.55,ymax=1.25,axis lines=left,
          xlabel={step $\tau$},ylabel={$\sigma_t A_t(P+\tau d)/\amin$},tick label style={font=\scriptsize},
          label style={font=\small},xtick={0,1},ytick={0.6,1}, samples=100, domain=0:1]
        \addplot[black!35,thin] {1 + 1.2*x - 1.6*x^2};
        \addplot[black!35,thin] {1 + 0.5*x - 0.1*x^2};
        \addplot[black!35,thin] {1.04 + 0.9*x - 2.4*x^2};
        \addplot[black!35,thin] {1.15 - 0.25*x};
        \addplot[black!35,thin] {1.02 + 0.3*x + 0.2*x^2};
        \addplot[accent,very thick,domain=0:1] {min(1 + 1.2*x - 1.6*x^2, 1 + 0.5*x - 0.1*x^2, 1.04 + 0.9*x - 2.4*x^2, 1.15 - 0.25*x, 1.02 + 0.3*x + 0.2*x^2)};
        \addplot[hot,only marks,mark=*,mark size=1.8pt] coordinates {(0.2189,1.0953)};
        \node[font=\scriptsize,anchor=south west,text=hot!80!black] at (axis cs:0.235,1.10) {step};
      \end{axis}
    \end{tikzpicture}
    \caption{Along an ascent direction each triangle area is a parabola in $\tau$; the step goes to the first local
    maximum of their lower envelope (blue). Schematic.}
  \end{subfigure}
  \caption{The two local moves on the nonsmooth objective: relocation of one point and the min-norm ascent step.}
  \label{fig:reloc}
\end{figure}

A configuration where no single point can be improved may still admit a joint move. Let $u$ be the free parameters
(Section~\ref{sec:setting}), $T_\varepsilon = \{t : a_t \le (1+\varepsilon)\amin\}$ and
$g_t = \nabla_u (\sigma_t A_t)$ (at most six non-zeros).
Following Demyanov and Malozemov \cite{DM74} we compute the minimum-norm element of the convex hull,
\[
  g^\ast = \arg\min \bigl\{ \|g\| : g \in \conv\{ g_t : t \in T_\varepsilon \} \bigr\},
  \qquad\text{which satisfies}\qquad \langle g^\ast, g_t \rangle \ge \|g^\ast\|^2 \quad \forall\, t \in T_\varepsilon
\]
by the optimality condition of the projection of $0$ onto a convex set. Hence, if $g^\ast \neq 0$, the direction
$d = g^\ast$ raises every $\varepsilon$-active triangle to first order; for $\varepsilon = 0$, $g^\ast = 0$ is exactly
Clarke stationarity of $\amin$ \cite{Clarke83}. We compute $g^\ast$ by pairwise Frank--Wolfe iterations over the
simplex of weights \cite{FW56,LJJ15}. The slack $\varepsilon > 0$ matters: active sets with $\varepsilon = 0$ jam at the kinks of
$\amin$. Along the ray $u + \tau d$ every $A_t$ is a quadratic polynomial in $\tau$, so the lower envelope
$\min_t \sigma_t A_t$ of parabolas is known in closed form, and the step goes to its first local maximum in
$\tau > 0$ (Figure~\ref{fig:reloc}). In practice relocation comes first and
ascent second: ascent from a raw start crawls, but it is the move that leaves a relocation-stationary point.

\subsection{A smooth surrogate}\label{sec:softmin}

For exploration (inside ruin and recreate and on the GPU) we use the soft minimum over a near-active subset
$S = \{t : a_t \le 6\amin(P_0)\}$, fixed at the start point $P_0$,
\begin{equation}\label{eq:softmin}
  F_\beta(P) \;=\; -\frac{\amin^S(P)}{c} + \frac1\beta \log \sum_{t \in S} \exp\Bigl(-\beta\,\frac{a_t(P) - \amin^S(P)}{c}\Bigr),
  \qquad c = \amin(P_0),
\end{equation}
with $\amin^S = \min_{t \in S} a_t$, which equals $-\operatorname{softmin}_\beta(a_t/c)$ and satisfies $-\amin^S/c \le F_\beta \le -\amin^S/c + \log|S|/\beta$.
In the convex variant $a_t$ is replaced by $a_t/H(P)$, whose gradient contains the term $-a_t \nabla H / H^2$, non-zero
on hull vertices only. $F_\beta$ is minimized by L-BFGS \cite{LBFGS} (memory $8$, Armijo backtracking) at
$\beta = 15, 20, 50$, with a quadratic penalty for leaving the container. The value $v$ of the resulting
configuration, computed in single precision, is its \emph{soft value}. Soft values lie a few percent below what the SLP
polish of Section~\ref{sec:slp} reaches from the same configuration; they rank starting points, they never decide a
record.

\subsection{Trust-region SLP with a lossless reduction}\label{sec:slp}

Every candidate passes the same polish. One step solves the LP
\begin{equation}\label{eq:slp}
  \max_{d,\,z}\ z \quad\text{s.t.}\quad v_t + \nabla v_t \cdot d \ge z \ \ (t \in R), \qquad \|d\|_\infty \le \delta,
  \qquad P + d \in K,
\end{equation}
over all $2n$ coordinates, with $v_t = \sigma_t A_t/\area(K)$ linearized at $P$ and a row set $R \subseteq T$ chosen
below. The step is accepted if the true value increases, and $\delta$ is adapted by the ratio of actual to predicted
gain (doubled above $3/4$ up to at most $0.1$, halved below $1/4$, multiplied by $0.3$ on rejection, halved when the LP
fails; $\delta_0 = 0.02/\sqrt n$).
The LPs are solved with HiGHS \cite{HiGHS}. For $n = 100$ there are $1.6 \cdot 10^5$ triangles, but very few matter
inside the trust region. Let $L_t$ be the sum of $|\Delta x| + |\Delta y|$ over the three edges of $t$. Since
$\partial A_t / \partial p_i$ is half the opposite edge rotated by $90^\circ$, $|\nabla A_t \cdot d| \le \frac12
\delta L_t$ whenever $\|d\|_\infty \le \delta$.

\begin{lemma}\label{lem:reduce}
In the square and the triangle let $q$ attain $v_q = v$. All rows $t$ with
\[
  v_t \;>\; v_q + \frac{\delta\,(L_t + L_q)}{2\,\area(K)}
\]
can be removed from \eqref{eq:slp} simultaneously without changing its feasible set, hence its optimal value and its set
of optimal solutions.
\end{lemma}

\begin{proof}
Row $q$ is never removed. Let $(d,z)$ satisfy row $q$ and the box. For every removed $t$,
\[
  z \;\le\; v_q + \frac{\delta L_q}{2\area(K)} \;<\; v_t - \frac{\delta L_t}{2\area(K)} \;\le\; v_t + \nabla v_t \cdot d,
\]
so the removed rows hold at every point that satisfies the kept rows.
\end{proof}

In the convex variant the same argument applies to the linearization of $a_t / H$, with
$|\nabla (a_t/H) \cdot d| \le \delta (L_t / (2H) + v_t \|\nabla H\|_1 / H)$. With this reduction the LPs have a few
hundred sparse rows (seven non-zeros each, plus the hull vertices in the convex variant), which makes the polish cheap
enough to be applied to every candidate the GPU forwards.

\subsection{The tie system and the exact gate}\label{sec:tie}

At a local maximum many triangles share the smallest value. Let $A \subseteq T$ be the active set and $u$ the free
parameters: two per interior point and one per side point in the square and the triangle; in the convex variant all
coordinates except those of the three points spanning the largest triangle, which are held fixed and thereby remove
the six-dimensional affine gauge (Lemma~\ref{lem:affine}), so $\dim u = 2n - 6$. We solve the \emph{tie system}
\begin{equation}\label{eq:tie}
  \frac{\sigma_t\, A_t(u)}{\area(K)} = \tau \qquad (t \in A), \qquad \area(K) = H(u) \text{ in the convex variant},
\end{equation}
in the unknowns $(u, \tau)$ by Gauss--Newton with the analytic Jacobian (minimum-norm steps when \eqref{eq:tie} is
underdetermined, least squares when it is overdetermined), in \texttt{mpmath} \cite{mpmath} with 57 significant digits.
The system is close to square: among the 184 stored configurations with $37 \le n \le 100$, the number of triangles
within relative distance $10^{-9}$ of the minimum equals the number of unknowns $\dim u + 1$ in 71 cases and differs
by at most three in 125. Because one threshold does not always identify $A$, we solve \eqref{eq:tie} for $A = T_\varepsilon$ with
$\varepsilon \in \{10^{-9}, 10^{-7}, 10^{-6}, 10^{-5}\}$, accept a solution only if no triangle outside $A$ ends below
$\tau$, and keep the best accepted $\tau$. The coordinates are then truncated toward zero to 32 decimals and the value
of exactly these decimals is recomputed in rational arithmetic, together with containment. This number is what the
record tables verify \cite{HSite}; a refinement that loses value against the double-precision polish is discarded. The
table of best configurations itself compares double-precision values and keeps a candidate only if it exceeds the
stored value by a relative $10^{-9}$; the tie system and the rational value are computed when a configuration is
exported for submission.

\section{Global moves}\label{sec:global}

\subsection{Stress-weighted ruin and recreate}\label{sec:rr}

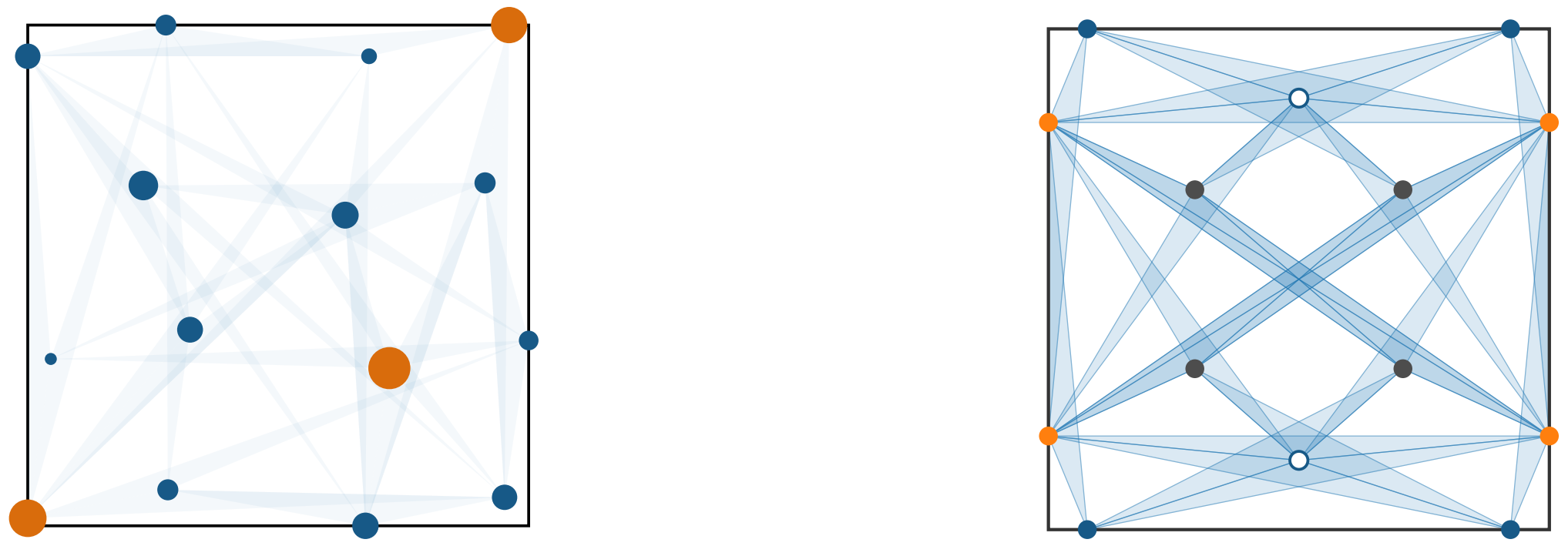
\begin{figure}[t]
  \centering
  \begin{subfigure}[b]{0.45\linewidth}\centering
    \begin{tikzpicture}[scale=4.3]
\draw[black,line width=0.7pt] (0,0) rectangle (1,1);
\fill[mplC0,fill opacity=0.05] (1.0000,0.3701) -- (0.9519,0.0569) -- (0.9129,0.6848) -- cycle;
\fill[mplC0,fill opacity=0.05] (1.0000,0.3701) -- (0.0000,0.0151) -- (0.2796,0.0719) -- cycle;
\fill[mplC0,fill opacity=0.05] (1.0000,0.3701) -- (0.0000,0.9377) -- (0.6337,0.6204) -- cycle;
\fill[mplC0,fill opacity=0.05] (1.0000,0.3701) -- (0.0460,0.3332) -- (0.7219,0.3149) -- cycle;
\fill[mplC0,fill opacity=0.05] (0.6741,0.0000) -- (0.9519,0.0569) -- (0.2796,0.0719) -- cycle;
\fill[mplC0,fill opacity=0.05] (0.6741,0.0000) -- (0.9609,1.0000) -- (0.9129,0.6848) -- cycle;
\fill[mplC0,fill opacity=0.05] (0.6741,0.0000) -- (0.0000,0.9377) -- (0.2308,0.6796) -- cycle;
\fill[mplC0,fill opacity=0.05] (0.6741,0.0000) -- (0.6337,0.6204) -- (0.7219,0.3149) -- cycle;
\fill[mplC0,fill opacity=0.05] (0.6741,0.0000) -- (0.6337,0.6204) -- (0.6814,0.9377) -- cycle;
\fill[mplC0,fill opacity=0.05] (0.6741,0.0000) -- (0.7219,0.3149) -- (0.9129,0.6848) -- cycle;
\fill[mplC0,fill opacity=0.05] (0.9519,0.0569) -- (0.9609,1.0000) -- (0.9129,0.6848) -- cycle;
\fill[mplC0,fill opacity=0.05] (0.9519,0.0569) -- (0.2757,1.0000) -- (0.7219,0.3149) -- cycle;
\fill[mplC0,fill opacity=0.05] (0.9519,0.0569) -- (0.0000,0.0151) -- (0.2796,0.0719) -- cycle;
\fill[mplC0,fill opacity=0.05] (0.9519,0.0569) -- (0.0000,0.9377) -- (0.2308,0.6796) -- cycle;
\fill[mplC0,fill opacity=0.05] (0.9609,1.0000) -- (0.0000,0.0151) -- (0.6337,0.6204) -- cycle;
\fill[mplC0,fill opacity=0.05] (0.9609,1.0000) -- (0.0000,0.9377) -- (0.6814,0.9377) -- cycle;
\fill[mplC0,fill opacity=0.05] (0.2757,1.0000) -- (0.0000,0.0151) -- (0.0460,0.3332) -- cycle;
\fill[mplC0,fill opacity=0.05] (0.2757,1.0000) -- (0.0000,0.9377) -- (0.6814,0.9377) -- cycle;
\fill[mplC0,fill opacity=0.05] (0.2757,1.0000) -- (0.3240,0.3916) -- (0.2796,0.0719) -- cycle;
\fill[mplC0,fill opacity=0.05] (0.0000,0.0151) -- (0.0000,0.9377) -- (0.0460,0.3332) -- cycle;
\fill[mplC0,fill opacity=0.05] (0.0000,0.0151) -- (0.6337,0.6204) -- (0.3240,0.3916) -- cycle;
\fill[mplC0,fill opacity=0.05] (0.0000,0.0151) -- (0.3240,0.3916) -- (0.6814,0.9377) -- cycle;
\fill[mplC0,fill opacity=0.05] (0.0000,0.9377) -- (0.3240,0.3916) -- (0.2308,0.6796) -- cycle;
\fill[mplC0,fill opacity=0.05] (0.6337,0.6204) -- (0.0460,0.3332) -- (0.9129,0.6848) -- cycle;
\fill[mplC0,fill opacity=0.05] (0.6337,0.6204) -- (0.2308,0.6796) -- (0.9129,0.6848) -- cycle;
\fill[accentdark] (1.0000,0.3701) circle (0.0196);
\fill[accentdark] (0.6741,0.0000) circle (0.0264);
\fill[accentdark] (0.9519,0.0569) circle (0.0253);
\fill[mplC1!85!black] (0.9609,1.0000) circle (0.0362);
\fill[accentdark] (0.2757,1.0000) circle (0.0208);
\fill[mplC1!85!black] (0.0000,0.0151) circle (0.0375);
\fill[accentdark] (0.0000,0.9377) circle (0.0254);
\fill[accentdark] (0.6337,0.6204) circle (0.0269);
\fill[accentdark] (0.0460,0.3332) circle (0.0120);
\fill[mplC1!85!black] (0.7219,0.3149) circle (0.0420);
\fill[accentdark] (0.3240,0.3916) circle (0.0258);
\fill[accentdark] (0.2308,0.6796) circle (0.0295);
\fill[accentdark] (0.2796,0.0719) circle (0.0210);
\fill[accentdark] (0.6814,0.9377) circle (0.0157);
\fill[accentdark] (0.9129,0.6848) circle (0.0211);\end{tikzpicture}
    \caption{Stress: square record $n=15$ with its 25 smallest triangles; disc size shows $s_i$, the three most
    stressed points in orange.}
  \end{subfigure}\hfill
  \begin{subfigure}[b]{0.45\linewidth}\centering
    \begin{tikzpicture}[scale=4.3]
\fill[cfgtri] (0.0776,0.0000) -- (0.0000,0.1869) -- (0.0000,0.8131) -- cycle;
\fill[cfgtri] (0.0776,0.0000) -- (1.0000,0.1869) -- (0.5000,0.1383) -- cycle;
\fill[cfgtri] (0.0776,0.0000) -- (0.7077,0.3213) -- (0.5000,0.1383) -- cycle;
\fill[cfgtri] (0.9224,1.0000) -- (1.0000,0.8131) -- (1.0000,0.1869) -- cycle;
\fill[cfgtri] (0.9224,1.0000) -- (0.0000,0.8131) -- (0.5000,0.8617) -- cycle;
\fill[cfgtri] (0.9224,1.0000) -- (0.2923,0.6787) -- (0.5000,0.8617) -- cycle;
\fill[cfgtri] (0.9224,0.0000) -- (0.0000,0.1869) -- (0.5000,0.1383) -- cycle;
\fill[cfgtri] (0.9224,0.0000) -- (1.0000,0.8131) -- (1.0000,0.1869) -- cycle;
\fill[cfgtri] (0.9224,0.0000) -- (0.2923,0.3213) -- (0.5000,0.1383) -- cycle;
\fill[cfgtri] (0.0776,1.0000) -- (0.0000,0.1869) -- (0.0000,0.8131) -- cycle;
\fill[cfgtri] (0.0776,1.0000) -- (1.0000,0.8131) -- (0.5000,0.8617) -- cycle;
\fill[cfgtri] (0.0776,1.0000) -- (0.7077,0.6787) -- (0.5000,0.8617) -- cycle;
\fill[cfgtri] (0.0000,0.1869) -- (1.0000,0.8131) -- (0.2923,0.3213) -- cycle;
\fill[cfgtri] (0.0000,0.1869) -- (1.0000,0.8131) -- (0.7077,0.6787) -- cycle;
\fill[cfgtri] (0.0000,0.1869) -- (1.0000,0.1869) -- (0.5000,0.1383) -- cycle;
\fill[cfgtri] (0.0000,0.1869) -- (0.2923,0.3213) -- (0.7077,0.6787) -- cycle;
\fill[cfgtri] (0.0000,0.1869) -- (0.2923,0.6787) -- (0.5000,0.8617) -- cycle;
\fill[cfgtri] (1.0000,0.8131) -- (0.0000,0.8131) -- (0.5000,0.8617) -- cycle;
\fill[cfgtri] (1.0000,0.8131) -- (0.2923,0.3213) -- (0.7077,0.6787) -- cycle;
\fill[cfgtri] (1.0000,0.8131) -- (0.7077,0.3213) -- (0.5000,0.1383) -- cycle;
\fill[cfgtri] (1.0000,0.1869) -- (0.0000,0.8131) -- (0.7077,0.3213) -- cycle;
\fill[cfgtri] (1.0000,0.1869) -- (0.0000,0.8131) -- (0.2923,0.6787) -- cycle;
\fill[cfgtri] (1.0000,0.1869) -- (0.7077,0.6787) -- (0.5000,0.8617) -- cycle;
\fill[cfgtri] (1.0000,0.1869) -- (0.7077,0.3213) -- (0.2923,0.6787) -- cycle;
\fill[cfgtri] (0.0000,0.8131) -- (0.2923,0.3213) -- (0.5000,0.1383) -- cycle;
\fill[cfgtri] (0.0000,0.8131) -- (0.7077,0.3213) -- (0.2923,0.6787) -- cycle;
\draw[cfgbox] (0,0) rectangle (1,1);
\node[orb0] at (0.0776,0.0000) {};
\node[orb0] at (0.9224,1.0000) {};
\node[orb0] at (0.9224,0.0000) {};
\node[orb0] at (0.0776,1.0000) {};
\node[orb1] at (0.0000,0.1869) {};
\node[orb1] at (1.0000,0.8131) {};
\node[orb1] at (1.0000,0.1869) {};
\node[orb1] at (0.0000,0.8131) {};
\node[orb2] at (0.2923,0.3213) {};
\node[orb2] at (0.7077,0.6787) {};
\node[orb2] at (0.7077,0.3213) {};
\node[orb2] at (0.2923,0.6787) {};
\node[orb3] at (0.5000,0.1383) {};
\node[orb3] at (0.5000,0.8617) {};\end{tikzpicture}
    \caption{Orbits: square record $n=14$, group $D_2$, orbits $4+4+4+2$ (one marker per orbit): five parameters
    instead of 20 (eight of the points lie on a side).}
  \end{subfigure}
  \caption{The two global devices: aimed ruin and symmetric starts.}
  \label{fig:global}
\end{figure}

Ruin and recreate \cite{RR00} alternates the destruction of part of a solution with its greedy reconstruction. For the
Heilbronn problem the destruction is aimed. Every triangle receives a Boltzmann weight and every point a stress,
\[
  w_t = \exp\Bigl(-\kappa\, \frac{a_t - \amin}{\amin}\Bigr), \qquad s_i = U_i \sum_{t \ni i} w_t, \qquad
  \kappa = 15, \ U_i \sim \mathcal{U}(0.9, 1.1),
\]
a soft count of the nearly smallest triangles through $p_i$ (Figure~\ref{fig:global}). A ruin removes the $k$ most
stressed points, $k = \lfloor n\, \mathcal{U}(0.10, 0.25) \rfloor$. The recreate step inserts points greedily up to
$n + o$ points, $o \in \{0, 1, 2\}$, polishing \eqref{eq:softmin} after each insertion, and then ruins back to $n$.
A move is accepted if it improves the soft value, or with probability $0.15$; the search restarts from the incumbent
every 30 seconds. A soft value above the incumbent is sharpened ($\beta = 8000$) and, in the square and the triangle,
refined by relocation sweeps and ascent (Sections~\ref{sec:reloc} and~\ref{sec:ascent}). The scheme follows the evolved program of
C.~Nemri \cite{cnemri} and the stress-weighted ruin introduced in the record tables by R.~Chappell \cite{HSite}.

\paragraph{Greedy insertion.} A new point is the best of $4000$ random candidates $q$ by
$\min_{j<k} |A_{(q,j,k)}|$, i.e.\ the candidate farthest, in the area sense, from all lines through two points.
Candidates are interior points, points on a side (probability $0.2$) or free corners ($0.002$). Storing the lines
$(\alpha_{jk}, \beta_{jk}, c_{jk})$ sorted by pair length lets a candidate be discarded as soon as it falls below the
best so far. In the convex variant candidates come from the bounding box enlarged by $5\%$, scored by
$\min_{j<k} |A_{(q,j,k)}| / H(P \cup \{q\})$, where $H$ grows by the triangles over the hull edges visible from $q$.

\subsection{Symmetric starts}\label{sec:sym}

Symmetric point sets have measure zero, so unconstrained multistart never finds them, yet many records are symmetric
(Figure~\ref{fig:gallery}). For cold starts the engine therefore works in orbit space: a point group $G$ compatible with
the region acts on orbit representatives, and the configuration is the union of their orbits ($n$ fixes the orbit
type counts). This divides the number of free parameters by about $|G|$. We use $C_2, C_4, D_1, D_2, D_4$ and diagonal
variants in the square ($D_4$ with vacancies or extra free points when $4 \nmid n$), $C_3, D_1, D_3$ in the triangle,
and $C_k$ ($k \le 6$) and $D_k$ ($k \le 12$) in the convex variant, sampled with weights from their observed yield. The best symmetric
outputs are then relaxed in the full space, and seeded symmetric runs move a symmetric configuration of the same group
from $n \pm |G|$ to $n$ by adding or removing an orbit.

\subsection{Moves along \texorpdfstring{$n$}{n} and between variants}\label{sec:moves}

Write $B_V(m)$ for the best configuration known to us for variant $V$ and $m$ points (our own, or a published record).
For the cell $(V, n)$ we use five kinds of starts:
\begin{center}
\begin{tabular}{ll}
  \textbf{up}: $B(n-1)$ plus one inserted point, & \textbf{up2}: $B(n-2)$ plus two, \\
  \textbf{dn}: $B(n+1)$ minus one point, & \textbf{dn2}: $B(n+2)$ minus two random points, \\
  \multicolumn{2}{l}{\textbf{mv}: $B(n)$ with one to three random points removed and re-inserted.}
\end{tabular}
\end{center}

\paragraph{Deletion scan.} For \textbf{dn} an exhaustive scan is affordable on the CPU: each of the $n+1$
deletions from $B(n+1)$ is polished (300 SLP iterations) and the best three are polished further (1500 iterations).
Ranking deletions without polishing does not work, since at a local maximum the smallest remaining area is nearly the
same for most deletions; after polishing only a few deletions gain substantially. The scan runs from the top cell down
and repeats whenever $B(n+1)$ has changed.

\paragraph{Transfer between variants.} By Lemma~\ref{lem:affine} a square or triangle configuration is a convex
configuration of at least the same value. Conversely, a convex configuration is mapped into the square by the affine
map that sends a smallest enclosing parallelogram of $\conv P$ onto $[0,1]^2$; such a parallelogram can be chosen with
a side of each pair of parallel sides containing an edge of the polygon \cite{STV95}, so it is found by trying all pairs
of edge directions. Into the triangle we map the smallest enclosing triangle bounded by three edge lines; since a
minimal enclosing triangle need only have one side flush with an edge \cite{ORourke86}, this is a heuristic. Before
the polish in the target variant the value drops by the factor $H(P)/\area(\text{enclosing set})$. For a nearly round
hull this factor is about $\pi/4 \approx 0.79$ for the square and $\pi/(3\sqrt3) \approx 0.60$ for the triangle; for
our convex configurations the medians are $0.80$ and $0.62$, a loss the polish has to recover.

\section{The GPU stage}\label{sec:gpu}

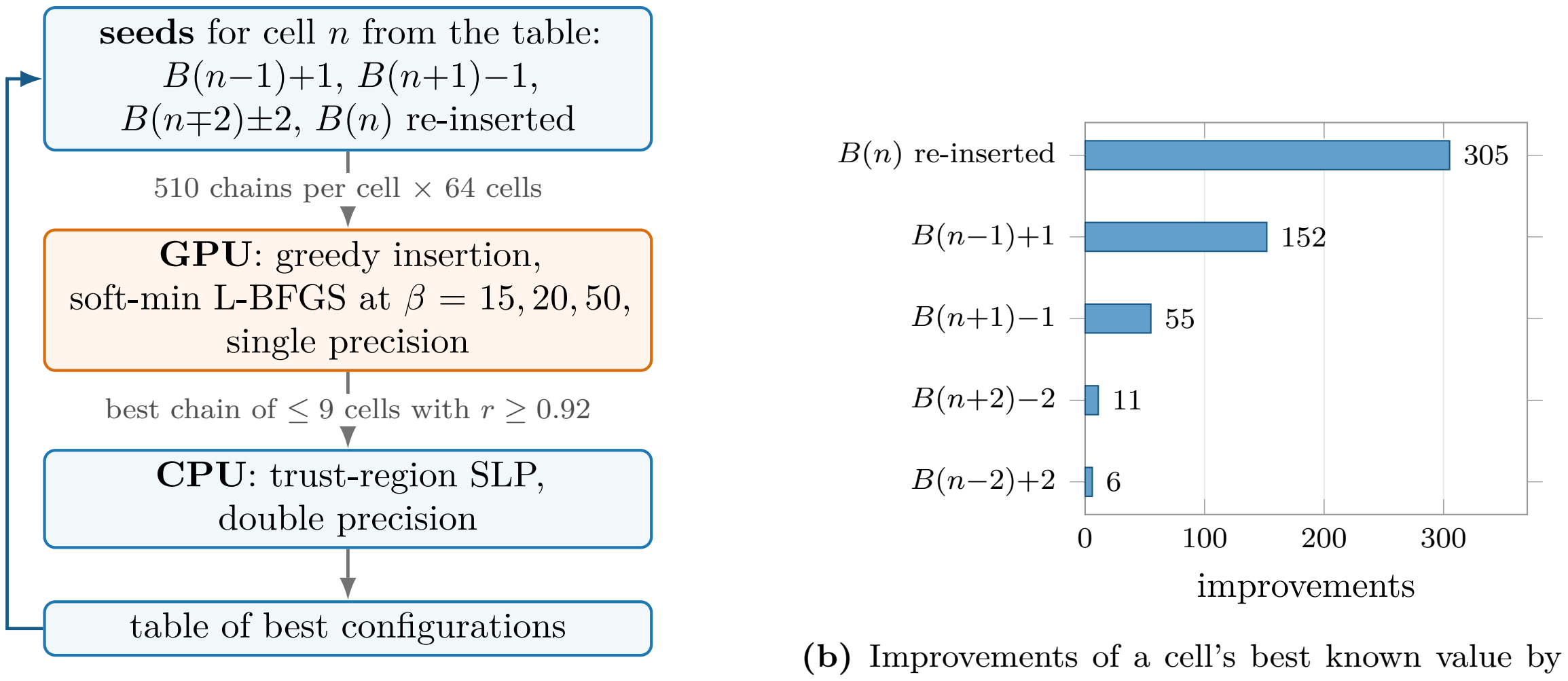
\begin{figure}[t]
  \centering
  \begin{subfigure}[b]{0.5\linewidth}\centering
    \begin{tikzpicture}[
        gbox/.style={rounded corners=3pt,draw=accent,fill=accent!6,line width=0.8pt,font=\small,text width=5.7cm,align=center,inner sep=3pt},
        ghot/.style={gbox,draw=hot!85!black,fill=hot!8},
        lab/.style={font=\scriptsize,text=black!70,fill=white,inner sep=1pt,align=center}]
      \node[gbox] (s) {\textbf{seeds} for cell $n$ from the table:\\ $B(n{-}1){+}1$, $B(n{+}1){-}1$, $B(n{\mp}2){\pm}2$, $B(n)$ re-inserted};
      \node[ghot,below=0.75cm of s] (g) {\textbf{GPU}: greedy insertion,\\ soft-min L-BFGS at $\beta = 15, 20, 50$,\\ single precision};
      \node[gbox,below=0.75cm of g] (c) {\textbf{CPU}: trust-region SLP,\\ double precision};
      \node[gbox,below=0.5cm of c] (r) {table of best configurations};
      \draw[flow] (s) -- node[lab] {$510$ chains per cell $\times$ $64$ cells} (g);
      \draw[flow] (g) -- node[lab] {best chain of $\le 9$ cells with $r \ge 0.92$} (c);
      \draw[flow] (c) -- (r);
      \draw[flow,draw=accentdark] (r.west) -- ++(-0.35,0) |- (s.west);
    \end{tikzpicture}
    \caption{One screening cycle.}
  \end{subfigure}\hfill
  \begin{subfigure}[b]{0.46\linewidth}\centering
    \begin{tikzpicture}
      \begin{axis}[xbar, width=0.8\linewidth, height=5.4cm, bar width=8pt, xmin=0, xmax=370,
          xlabel={improvements}, symbolic y coords={up2,dn2,dn,up,mv}, ytick=data,
          yticklabels={{$B(n{-}2){+}2$},{$B(n{+}2){-}2$},{$B(n{+}1){-}1$},{$B(n{-}1){+}1$},{$B(n)$ re-inserted}},
          nodes near coords, nodes near coords style={font=\scriptsize}, tick label style={font=\scriptsize},
          label style={font=\small}, xmajorgrids, grid style={gray!15}, axis line style={black!40}]
        \addplot[fill=accent!70,draw=accentdark] coordinates {(6,up2) (11,dn2) (55,dn) (152,up) (305,mv)};
      \end{axis}
    \end{tikzpicture}
    \caption{Improvements of a cell's best known value by seed type, 26~Nov.\ 09:20 to 27~Nov.\ 11:10.}
  \end{subfigure}
  \caption{GPU screening: cheap in, expensive out.}
  \label{fig:gpu}
\end{figure}

The starts of Section~\ref{sec:moves} are cheap but mostly useless. The GPU stage (Figure~\ref{fig:gpu}) screens them
in bulk; it is the moves along $n$ of Section~\ref{sec:moves} at scale.

\paragraph{Batch.} For every cell $n = 37, \dots, 100$ of a variant, $510$ chains are seeded per cycle, equally split
over \textbf{up}, \textbf{dn}, \textbf{up2}, \textbf{dn2}, \textbf{mv}; for 64 cells these are $32\,640$ chains, one GPU
thread each, sorted by $n$ so that SIMD groups stay uniform. Each chain is completed to $n$ points by greedy insertion
(each insertion followed by $20$ L-BFGS iterations of \eqref{eq:softmin} at $\beta = 15$), then polished at $\beta = 20$
and $\beta = 50$ ($50$ iterations each), all in single precision. The kernels are written in Metal; the convex kernel
adds the hull area and its gradient (monotone chain, with the $x$-order of the points kept between calls).

\paragraph{Gate.} For each cell let $r$ be the best soft value of its chains divided by the current target, the larger
of the stored and the published value. Cells with $r \ge 0.92$ are candidates; at most nine per cycle are forwarded,
ordered first by the number of SLP polishes since the cell's last real gain (above a relative $10^{-5}$; this keeps a
cell that produces near misses from monopolizing the CPU) and then by $r$. The best chain of each forwarded cell
receives the SLP polish (300 iterations, double precision), and a result above the stored value enters the table. A cycle over the 64 cells of one variant takes 20--35 seconds.

\paragraph{What did not help.} A first port that ran the complete ruin and recreate search on the GPU, one chain per
thread, produced two improvements in eleven hours. Further ruin and recreate steps after the seeding did not raise the
yield, and ranking deletions by their soft value instead of polishing each of them failed. The working design
uses the GPU for what it is good at, evaluating very many independent candidates, and leaves every decision on a value
to double precision and, finally, to rational arithmetic. In the 26 hours after its start (26 November 2026, all three
variants, the convex kernel from the second morning), 529 of the 643 improvements of a best known value came from GPU
seeds, 81 from deletion scans and 33 from the follow-up seeds of earlier improvements.

\section{Implementation}\label{sec:impl}

The released code \cite{Code} is a compact version of the research code that produced the records; Table~\ref{tab:code}
maps the sections to its modules. A Python driver keeps the table of best configurations (one file per cell) and runs
three loops concurrently, all reading their starts from the table: GPU screening cycles over the variants, deletion
scans from the top cell down, and C engine runs (cold symmetric and warm ruin and recreate), whose outputs are polished
and offered to the table. Every improvement is also transferred to the other variants. (Seeded symmetric runs and
hot basin-hopping walks are options of the engine that the driver does not schedule.) CPU work runs in one process
pool, so the machine is never oversubscribed; the records were found on a single workstation (Apple M5 Pro, 15 CPU
cores, integrated GPU). Without a GPU the driver replaces the screening by greedy grid insertion on the CPU.

\begin{table}[t]
  \centering\small
  \begin{tabular}{lll}
    \toprule
    component & section & module \\
    \midrule
    one-point LP relocation, min-norm ascent & \ref{sec:reloc}, \ref{sec:ascent} & C: \texttt{reloc.c} \\
    soft-min L-BFGS, ruin and recreate & \ref{sec:softmin}, \ref{sec:rr} & C: engine \texttt{heil} \\
    symmetric starts & \ref{sec:sym} & C: engine \texttt{heil} \\
    trust-region SLP with lossless reduction & \ref{sec:slp} & \texttt{polish.py} \\
    tie system, truncation, rational value & \ref{sec:tie} & \texttt{refine.py}, \texttt{exact.py} \\
    moves along $n$, deletion scan & \ref{sec:moves} & \texttt{neighbours.py} \\
    transfer between variants & \ref{sec:moves} & \texttt{transfer.py} \\
    GPU kernels (Metal) & \ref{sec:gpu} & \texttt{gpu.py} \\
    table of best configurations, driver & \ref{sec:impl} & \texttt{store.py}, \texttt{run.py} \\
    \bottomrule
  \end{tabular}
  \caption{Where the components live in the released code: C sources in \texttt{csrc/}, Python modules in \texttt{heilbronn/}.}
  \label{tab:code}
\end{table}

\section{Results}\label{sec:results}

On \StateDate{} our configurations were the records in \NOurs{} of the \NCells{} cells with $3 \le n \le 100$ of the
public tables \cite{HSite}, and in \NOursBig{} of the \NBig{} cells with $37 \le n \le 100$: all $64$ triangle and all
$64$ convex cells, and \NOursBigSquare{} of the $64$ square cells. All configurations are available with exact coordinates on the
record site. Since the tables change daily, we do not list individual values.

Two observations are worth recording. First, the improvements travel through the table in waves: a new configuration
for $n$ changes the starts for $n \pm 1$ and $n \pm 2$ and for the other variants, and often several neighbouring cells
improve within minutes. The largest single jumps came from symmetric starts in the convex variant (for example, $D_{11}$
raised our value for $n = 100$ by $24.6\%$, and the seeds built from it raised $n = 99$ and $n = 98$ by about $20\%$).
Second, for small $n$ none of the moves helped during our campaign (November 2026): for $n \le 40$ neither seeding nor a
GPU random multistart with about $2500$ SLP polishes improved any record, and the deletion scan reproduced many known
configurations exactly. This is consistent with the records for small $n$ being optimal or very close to it, as is
certified computationally for $n \le 10$ in the square and $n \le 8$ in the triangle \cite{SM26a,SM26b}, whereas for $n \ge 60$ single starts still
gained several percent until recently.

\section{Discussion}\label{sec:discussion}

The most effective idea was to treat the table as a whole: good configurations for $n$ points carry most of the
structure of good configurations for $n \pm 1$, and a careful local polish recovers the rest. The second is the division
of labour between imprecise but massively parallel screening and local optimization of the nonsmooth objective in
double precision, whose LPs can be reduced without loss. Neither idea is specific to the Heilbronn problem; we expect
both to carry over to other max--min point-placement problems, such as the higher-dimensional variants.

The pipeline was built and run in close collaboration with an AI coding assistant. Our experience is that this works
well for a problem with a cheap, unambiguous and exactly checkable objective: the assistant can run, measure and
discard many ideas, while every claimed record is checked by an independent verifier in rational arithmetic.

\paragraph{Code and data.} Code: \CodeURL{}. Configurations: \url{https://math.tejstead.com/heilbronn}.


\end{document}